\documentclass[11pt]{article}
\usepackage{epic,latexsym,amssymb,xcolor}
\usepackage{color}
\usepackage{tikz}
\usetikzlibrary{shapes.geometric}
\usepackage{amsfonts,epsf,amsmath,leftidx}
\usepackage{float}
\usepackage{pgfplots}
\usepackage{mathrsfs}
\usetikzlibrary{arrows}
\usepackage{pdflscape}
\usepackage[bottom=1.5in, top=1.2in, left=1.5in, right=1.5in]{geometry}

\newtheorem{thm}{Theorem}

\newcommand{\qed}{$\Box$}

\let\oldenumerate\enumerate
\renewcommand{\enumerate}{
  \oldenumerate
  \setlength{\itemsep}{0pt}
  \setlength{\parskip}{0pt}
  \setlength{\parsep}{0pt}
}

\begin{document}

\title{A simple proof on the number of $(3 \times n)$-Latin rectangles based on a set of $\lambda$ elements}

\author{Pantaree Thengarnanchai$^{a,b,}$\thanks{Petchra Pra Jom Klao Master Scholarship from King Mongkut's University of Technology Thonburi~(10/2566)}, Pawaton Kaemawichanurat$^{a, b,}$\thanks{National Research Council of Thailand (NRCT) and King Mongkut's University of Technology Thonburi (N42A660926)},\\ Watcharintorn Ruksasakchai$^{c}$, Natawat Klamsakul$^{a, b}$
\\ \\
$^a$Department of Mathematics, Faculty of Science,\\
King Mongkut's University of Technology Thonburi, \\
Bangkok, Thailand \\
$^b$Mathematics and Statistics with Applications (MaSA) \\
$^c$Department of Computational Science and Digital Technology,\\
Faculty of Liberal Arts and Science,\\
Kasetsart University, Kamphaeng Saen Campus\\
\small \tt pantaree.theng@gmail.com, pawaton.kae@kmutt.ac.th,\\
\tt watcharintorn1@hotmail.com, natawat.kla@kmutt.ac.th}

\date{}
\maketitle

\begin{abstract}
In 1980, Athreya, Pranesachar and Singhi established the chromatic polynomial of $(3 \times n)$-Latin rectangles whose entries based on a set $\{1, 2, ..., \lambda\}$ in which $\lambda \geq n$. Their proof requires M\"{o}bius inversion formula and lattice partitions. In this paper, we present a simpler proof by using the idea of mathematical induction and appropriate coloring. 
\end{abstract}


\section{Introduction and Motivation}
Throughout this paper, we let a graph $G = (V(G), E(G))$ be finite and simple, no loops or multiple edges. For a vertex $x \in V(G)$, a \emph{neighbour} of $x$ in $G$ is a vertex $y$ such that $xy \in E(G)$. The \emph{neighbour set} of $x$ in $G$ is the set of all neighbours of $x$ in $G$ and is denoted by $N_{G}(x)$. For a graph $G$, the \emph{line graph} $H_{G}$ of $G$ is the graph whose vertices correspond to the edges of $G$ and any two vertices of $H$ are adjacent if and only if the corresponding edges of $G$ are incident. Let $u, v$ be a pair of distinct vertices of $G$. The \emph{identified graph} $G_{uv}$ is obtained by removing the vertices $u$ and $v$ and adding the new vertex $x_{uv}$,  and join $x_{uv}$ to all vertices in $(N_{G}(u) \cup N_{G}(v)) \setminus \{u, v\}$. For graphs $H$ and $G$, the \emph{lexicographic product} between $G$ and $H$ is the graph whose vertex set is the cartesian product $V(G) \times V(H)$ in which two vertices $(u, v)$ and $(u', v')$ are adjacent if (i) $u = u'$ and $vv' \in E(H)$ or (ii) $uu' \in E(G)$ and $v = v'$. For a given graph $G$ and a natural number $\lambda$, the \emph{chromatic polynomial} of $G$, denoted by $P(G, \lambda)$, is the function which gives the number of proper coloring of $G$ by using $\lambda$ colors. For a coloring function $f: V(G) \rightarrow \mathbb{N}$, we denote by $col(v) = i$ if a vertex $v \in V(G)$ is given the color $i \in \mathbb{N}$ by $f$.
\vskip 5 pt

\indent For natural numbers $n \geq m$, a \emph{Latin rectangle} is an array of dimension $m \times n$ whose entries are $1, 2, ..., n$, each of which appears exactly once in each row and at most once in each column. The problem of interest is to count all the possible Latin rectangles when one dimension of the array $m \times n$ is fixed. We let $\mathcal{L}_{m, n}$ be the family of Latin rectangles whose entries of the first row are $1, ..., n$, respectively. We let $L_{m, n} = |\mathcal{L}_{m, n}|$. Interestingly, $L_{2, n} = n!\sum_{k=0}^{n}\frac{(-1)^{k}}{k!}$, which is the derangement formula of $n$ objects. When $m = 3$, Riordan \cite{A} proved in 1944 that:
\begin{align}\label{l3n}
    L_{3, n} = (n!) \sum_{k+j \leq n}\frac{2^{j}}{j!}k!{-3(k+1) \choose n-k-j}.
\end{align}
\noindent Further, Arthreya et. al. \cite{APS} applied the concept of M\"{o}bius inversion formula and Lattice partition to establish the chromatic polynomial $g(n, \lambda)$ of $(3 \times n)$-Latin rectangle whose entries are members in the set $\{1, 2, ..., \lambda\}$ in which $\lambda \geq n$, each of the entries appears at most once in each column and row. They proved that:
\begin{align}\label{l3n}
g(n, \lambda) = \frac{\lambda!n!}{((\lambda - n)!)^{3}}\sum_{\alpha + \beta + \gamma = n}(-1)^{\beta}2^{\gamma}\frac{((\lambda - n + \alpha)!)^{2}}{\alpha!\gamma!}\binom{3\lambda - 3n + 3\alpha + \beta + 2}{\beta}.
\end{align}
\noindent In this paper, we present a simpler proof to establish $g(n, \lambda)$, requiring only the ideas of mathematics induction and appropriate coloring. Our proof might last several pages as it is clearly explained.


\section{The chromatic polynomials of the line graphs of $K_{3, n}$ with $\lambda$ colors}

\noindent We begin with a well-known tool, the so called \emph{reduction formula}, which is employed to establish the chromatic polynomial of graph.
\vskip 5 pt

\begin{thm}
Let $G$ be a graph with $uv \in E(G)$ and $\lambda$ be a natural number. We have that
\begin{align*}
 P(G, \lambda) = P(G - uv, \lambda) - P(G_{uv}, \lambda).   
\end{align*}
\end{thm}




\indent First of all we let $G(n)$ be the line graph of complete bipartite $K_{3,n}$. Thus $G(n) = K_{3} \Box K_{n}$ consists of 3 cliques $K_{n}^{1} , K_{n}^{2} , K_{n}^{3}$ of order 3. We let $K_{n}^{i} = \{ x_{1}^{i} , x_{2}^{i} , x_{3}^{i} , ... , x_{n}^{i} \}$. Since $G(n) = K_{3} \Box K_{n}$, we let without loss of generality that the set $\{ x_{j}^{1} , x_{j}^{2} , x_{j}^{3} \}$ form a triangle for each $1 \leq j \leq n$. 

\begin{figure}[h]
\centering
{\includegraphics[width=0.50\textwidth]{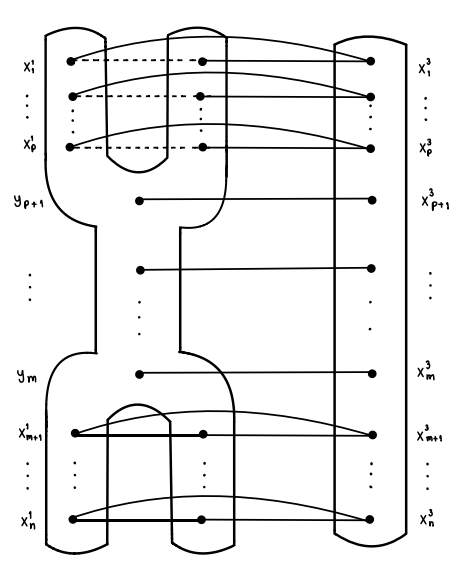}}%
\caption{\footnotesize{The graph $G(n, p, q)$.}}
\label{fig1}
\end{figure}

We construct the following graph from $G(n)$ for non-negative integers $n, m, p$ and $q$ such that $n \geq m = p + q \geq 0$, we let the graph $G(n,p,q)$ be obtained from $G(n)$ by removing the edges $x_{1}^{1}x_{1}^{2}, x_{2}^{1}x_{2}^{2} , ... , x_{p}^{1}x_{p}^{2}$ and identifying the edges $x_{p+1}^{1}x_{p+1}^{2} , ... , x_{m}^{1}x_{m}^{2}$ with the new vertices $y_{p+1} , ... , y_{m}$ respectively. The graph $G(n, p, q)$ is illustrated in Figure \ref{fig1}.
\vskip 15 pt

\vskip 15 pt

\begin{thm}\label{binomchro}
For natural numbers $\lambda$ and $n$ such that $\lambda \geq n$, we let $g(n,\lambda)$ be the chromatic polynomial of $G(n)$. Further, for non negative integers $m, p, q$ such that $n \geq m = p + q \geq 1$, we let $g(n,p,q,\lambda)$ be the chromatic polynomial of coloring the graph $G(n,p,q)$ with $\lambda$ colors. Then
\begin{align}\label{thm}
g(n,\lambda) = \sum_{q=0}^{m}{m \choose q}(-1)^{q}g(n,p,q,\lambda).
\end{align}
\end{thm}

\noindent \emph{Proof.} We will prove by induction on $m$ for all $1 \leq m \leq n$. When $m=1$, Equation (\ref{thm}) is proved by the reduction formula, establishing basic step. Hence, we may assume that Equation (\ref{thm}) holds for the case $m'$ when $m'=m - 1$. Therefore,
\begin{align}\label{induct}
g(n,\lambda) = \sum_{q'=0}^{m-1}{m-1 \choose q'}(-1)^{q'}g(n,p',q',\lambda)
\end{align}
\noindent where $p' = m - 1 - q'$. That is, for each $0 \leq q' \leq m-1$, the graph $G(n,p',q')$ has the chromatic polynomial $g(n,p',q', \lambda)$ with the coefficient ${m-1 \choose q'}(-1)^{q'}$

\indent By applying the reduction formula to $G(n,p',q')$ for each $0 \leq q' \leq m-1$ with the edge $x^{1}_{m}x^{2}_{m}$ removing and identifying, we have
\begin{align}\label{q=0}
g(n, m-1, 0, \lambda) = g(n,m,0,\lambda) - g(n,m - 1,1, \lambda)
\end{align}
when $q'=0$ and we have
\begin{align}\label{q>0}
g(n, p', q', \lambda) = g(n,p'+1,q',\lambda) - g(n,p',q'+1, \lambda)
\end{align}
when $1 \leq q' \leq m - 1$. By plugging Equations (\ref{q=0}) and (\ref{q>0}) to (\ref{induct}), we have that
\begin{align*}
g(n,\lambda) = &{m - 1 \choose 0}(-1)^{0}g(n,m - 1,0,\lambda) + {m - 1 \choose 1}(-1)^{1}g(n,m - 2,1,\lambda) \cdots \\
& + {m - 1 \choose m - 1}(-1)^{m - 1}g(n,0,m - 1,\lambda)\\
= &{m - 1 \choose 0}(-1)^{0}(g(n,m,0,\lambda) - g(n,m - 1,1, \lambda))\\
&+ {m - 1 \choose 1}(-1)^{1}(g(n,m - 1,1,\lambda) - g(n,m - 2,2,\lambda))\\
&+ {m - 1 \choose 2}(-1)^{2}(g(n,m - 2,2,\lambda) - g(n,m - 3,3,\lambda))\\
& \vdots\\
&+ {m - 1 \choose m-1}(-1)^{m-1}(g(n,1,m-1,\lambda) - g(n,0,m,\lambda))\\
= &{m - 1 \choose 0}(-1)^{0}g(n,m,0,\lambda)\\
&+ ({m -1 \choose 0} + {m - 1 \choose 1})(-1)^{1}g(n,m - 1,1,\lambda)\\
&+ ({m - 1 \choose 1} + {m - 1 \choose 2})(-1)^{2}g(n,m - 2,2,\lambda)\\
&+ ({m - 1 \choose 2} + {m - 1 \choose 3})(-1)^{3}g(n,m - 3,3,\lambda)\\
& \vdots\\
&+ {m - 1 \choose m-1}(-1)^{m}g(n,0,m,\lambda)
\end{align*}
\noindent which implies that
\begin{align*}
g(n,\lambda) = \sum_{q=0}^{m}{m \choose q}(-1)^{q}g(n,p,q,\lambda)
\end{align*}
because $\binom{m - 1}{0} = \binom{m}{0}, \binom{m - 1}{m - 1} = \binom{m}{m}$ and $\binom{r - 1}{s - 1} + \binom{r - 1}{s} = \binom{r}{s}$ for all natural numbers $r, s$. This proves Equation (\ref{thm}).
\qed
\indent

\indent Now, we are ready to establish the chromatic polynomial of $(3 \times n)$-Latin rectangles which the entries are $1, 2, ..., \lambda$ so that $\lambda \geq n$. In the following, we let $G(n, k, l)$ be the graph $G(n, p, q)$ such that $k = p, l = q$ but $k + l = n$. Further, we may define some generalization of derangement. For a non-negative integer $t$ such that $0 \leq t \leq n$, we let $D_{\lambda, n, t}$ to denote the number of derangement of elements in a given set $S \subseteq \{1, 2, ..., \lambda\}$ to the positions $1, 2, ..., n$ such that $|S \cap \{1, 2, ..., n\}| = t$. By the inclusion-exclusion principle, it can be showed that
\begin{align}\label{derangement}
D_{\lambda, n, t} = \sum_{i = 0}^{t}(-1)^{i}\frac{(\lambda - i)!}{(\lambda - n)!}\binom{t}{i}.    
\end{align}
\noindent We prove that:
\vskip 5 pt

\begin{thm}
For a natural numbers $\lambda \geq n \geq 1$, we let $g(n,\lambda)$ be the chromatic polynomial of $G(n)$, the line of $K_{3, n}$. Then
\begin{align}
g(n,\lambda)= \frac{\lambda!}{(\lambda-n)!}\sum_{l=0}^{n}(-1)^{l}{n \choose l} \sum_{b}^{t_{1} = 0}\sum_{l - t_{1}}^{t_{2} = 0}( A(\lambda, k, l, t_{1}, t_{2})(B(\lambda, k, l, t_{1}))^{2})
\end{align}
where $k = n - l$ and
\begin{align*}
A(\lambda, k, l, t_{1}, t_{2}) &= \binom{k}{t_{1}}\binom{l}{t_{2}}\binom{\lambda - n}{l - t_{1} - t_{2}}\sum^{t_{2}}_{i = 0}(-1)^{i}(l - i)!\binom{t_{2}}{i},\\
B(\lambda, k, l, t_{1}) &= \sum_{t_{3} = 0}^{k - t_{1}}(\binom{k - t_{1}}{t_{3}}\binom{\lambda - n + t_{1}}{k - t_{3}}\sum^{t_{3}}_{j = 0}(-1)^{j}(k - j)!\binom{t_{3}}{j}).
\end{align*}
\end{thm}

\noindent \emph{Proof.} By applying Theorem \ref{binomchro} when $k= p, l = q$ and $m = n$, it suffices to show that  
$$g(n, k, l, \lambda) = \frac{\lambda!}{(\lambda-n)!}\sum_{b}^{t_{1} = 0}\sum_{l - t_{1}}^{t_{2} = 0}( A(\lambda, k, l, t_{1}, t_{2})(B(\lambda, k, l, t_{1}))^{2})$$
for all $k + l = n$.
\vskip 5 pt

\indent For the sake of convenient, we partition the graph $G(n,k,l)$ as follows : \vskip 5 pt

\begin{itemize}   
   \item $A = \{ x_{1}^{3}, x_{2}^{3},..., x_{k}^{3} \}$
   \item $B = \{ x_{k+1}^{3}, x_{k+2}^{3},..., x_{n}^{3} \}$
   \item $C = \{ y_{k+1}, y_{k+2},..., y_{n} \}$ 
   \item $D = \{ x_{1}^{1}, x_{2}^{1},..., x_{k}^{1} \}$
   \item $E = \{ x_{1}^{2}, x_{2}^{2},..., x_{k}^{2} \}$
\end{itemize}

\begin{figure}[h]
\centering
{\includegraphics[width=0.30\textwidth]{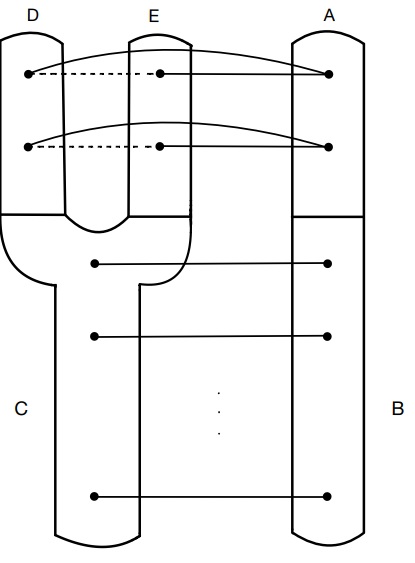}}%
\caption{\footnotesize{The partition of graph $G(n, k, l)$.}}
\label{fig2}
\end{figure}

\noindent The graph partition is illustrated in Figure \ref{fig2}. 
\vskip 15 pt

\vskip 15 pt

\noindent It can be checked that 
\begin{center}
 $\vert A \vert = \vert D \vert = \vert E \vert = k$ and $\vert B\vert = \vert C \vert = l$ and $k+l=n$   
\end{center}
\vskip 5 pt

\indent Let $\mathcal{C} = \{1, 2, ..., \lambda \}$ be the set of $\lambda$ colors that we will assign to the vertices of $G(n,k,l)$. We first color the vertices in $A \cup B$. There are $\frac{\lambda!}{(\lambda-n)!}$ possibilities. We next consider coloring on $C$ and then $D \cup E$. Without loss of generality we assume that 

\begin{center}
\noindent $col(x_{k+1}^{3})=1, col(x_{k+2}^{3})=2, ...,col(x_{n}^{3})=l$ and \vskip 5 pt
\noindent $col(x_{1}^{3})=l+1, col(x_{2}^{3})=l+2, ...,col(x_{k}^{3})=l+k$.
\end{center}
\vskip 5 pt

\indent Let $b = \min \{k, l\}$. For $0 \leq t_{1} \leq b$, we let $T_{1}$ be a subset of $\{l+1, l+2, ...,l+k\}$ of $t_{1}$ colors that are given to vertices in $C$ and, for $0 \leq t_{2} \leq l-t_{1}$, we let $T_{2}$ be a subset of $\{ 1, 2, ..., l \}$ of $t_{2}$ colors that are given to vertices in $C$. Next, we let $S$ be a subset of $ \{ l+k+1, l+k+2, ..., \lambda \}$ of $l-t_{1}-t_{2}$ colors that are given to vertices in $C$. There are ${k \choose t_{1}}{l \choose t_{2}}{\lambda -l-k \choose l-t_{1}-t_{2}}$ possibilities of choosing $T_{1}, T_{2}$ and $S$. Further, there are $D_{l,l,t_{2}}$ possibilities to color the vertices in $C$ by using the colors in $T_{1} \cup T_{2} \cup S$. Therefore, by Equation (\ref{derangement}), we can color the vertices in $C$ by
\begin{align*}
A(\lambda, k, l, t_{1}, t_{2}) &= \binom{k}{t_{1}}\binom{l}{t_{2}}\binom{\lambda - n}{l - t_{1} - t_{2}}D_{l, l, t_{2}}\\
&= \binom{k}{t_{1}}\binom{l}{t_{2}}\binom{\lambda - n}{l - t_{1} - t_{2}}\sum^{t_{2}}_{i = 0}(-1)^{i}(l - i)!\binom{t_{2}}{i}  
\end{align*}
\noindent possibilities.
\vskip 5 pt

\indent It can be checked that we can color independently between $D$ and $E$, each of which depending on the derangement of some color in $\{ l+1, l+2, ..., l+k \}$. For $0 \leq t_{3} \leq k-t_{1}$, we let $T_{3}$ be a subset of $\{ l+1,l+2, ..., l+k \}\setminus T_{1}$ of $t_{3}$ colors that are given to vertices in $D$ and $U$ be a subset of $\{ 1,2,...,\lambda \}\setminus (\{l+1, l+2,...l+k\}\cup T_{2} \cup S)$ of $k-t_{3}$ colors that are used to color vertices in $D$. Further, we can color the vertices in $D$ with the elements in $T_{3} \cup U$ in $D_{k, k, t_{3}}$ ways. Hence, we can color $D$ in ${k-t_{1} \choose t_{3}}{\lambda -k-(l-t_{1})\choose k-t_{3}}D_{k, k, t_{3}}$ ways. Similarly, we can color $E$ in ${k-t_{1} \choose t_{3}}{\lambda -k-(l-t_{1})\choose k-t_{3}}D_{k, k, t_{3}}$ ways. Hence, by Equation (\ref{derangement}), we can color the vertices in each of $D$ and $E$ by
\begin{align*}
B(\lambda, k, l, t_{1}) &= \sum_{t_{3} = 0}^{k - t_{1}}(\binom{k - t_{1}}{t_{3}}\binom{\lambda - n + t_{1}}{k - t_{3}}D_{k, k, t_{3}})\\
&= \sum_{t_{3} = 0}^{k - t_{1}}(\binom{k - t_{1}}{t_{3}}\binom{\lambda - n + t_{1}}{k - t_{3}}\sum^{t_{3}}_{j = 0}(-1)^{j}(k - j)!\binom{t_{3}}{j})    
\end{align*}
possibilities.
\vskip 5 pt

\noindent Therefore, we have 
\begin{align}
g(n,k,l)= \frac{\lambda!}{(\lambda-n)!}\sum_{l=0}^{n}(-1)^{l}{n \choose l} \sum_{b}^{t_{1} = 0}\sum_{l - t_{1}}^{t_{2} = 0}( A(\lambda, k, l, t_{1}, t_{2})(B(\lambda, k, l, t_{1}))^{2}).
\end{align}
This completes the proof.\qed


\noindent \textbf{Acknowledgment:} We would like to express our sincere thank to Ian Wanless from Monash University for his valuable suggestion, using the idea of chromatic polynomial to establish the formula of $g(n, \lambda)$.

\end{document}